\documentclass[11pt,a4paper]{article}
\usepackage{amsmath,amssymb,amsfonts}
\usepackage{a4wide}
\usepackage{epsfig,epic,eepic,graphicx}
\newtheorem{theorem}{Theorem}

\newtheorem{proposition}{Proposition}

\newcommand\bftau{\mbox{\boldmath$\tau$\unboldmath}}

\makeatletter

\@addtoreset{equation}{section}
\makeatother

\newcommand{\T}{{\mathbb T}}
\newcommand{\N}{{\mathbb N}}
\newcommand{\Z}{{\mathbb Z}}

\newcommand{\R}{{\mathbb R}}
\newcommand{\C}{{\mathbb C}}
\newcommand{\pa}{{\partial}}
\newcommand{\na}{{\nabla}}

\newcommand{\eps}{{\varepsilon}}

\def\Im{{\cal I}m \,}

\def\Re{{\rm Re}  \,}
\def\u{\mathbf{u}}
\def\x{\mathbf{x}}

\title{On the ill-posedness of the  Prandtl equation}
\author{\footnote{DMA/CNRS, Ecole Normale Sup\'erieure,
 45 rue d'Ulm,75005 Paris} David G\'erard-Varet, 
 \footnote{ENS/IPGP/CNRS,  Ecole Normale Sup\'erieure, 29 rue Lhomond,75005
   Paris} Emmanuel Dormy}

\date{}

\begin{document}
\maketitle

\begin{abstract}
The concern of this paper is the Cauchy problem for the Prandtl equation. 
This problem is known to be well-posed  for analytic data \cite{Sam:1998a,Lom:2003},
or for  data with monotonicity properties \cite{Ole:1999, Xin:2004}. We prove here that  
 it is linearly ill-posed in Sobolev type spaces. The key of the analysis is the
 construction, at high tangential frequencies, of unstable quasimodes for
 the linearization  around solutions with non-degenerate critical
 points. Interestingly, the strong instability is due to vicosity, which is coherent with
 well-posedness results obtained for the inviscid version of the
 equation \cite{Hon:2003}. A numerical study  of this  instability is also 
 provided.
\end{abstract}

\section{Introduction}
One  challenging open problem of fluid dynamics is to understand the
inviscid limit of the Navier-Stokes equations
\begin{equation} \label{NS}
\left\{
\begin{aligned}
\pa_t \u^\nu  + \u^\nu \cdot \na \u^\nu + \na p^\nu - \nu \Delta \u^\nu \:
 = \: 0, & 
\quad  \x \in \Omega, \\
\na \cdot \u^\nu \:  = \: 0, & \quad \x \in \Omega,\\
\u^\nu\vert_{\pa \Omega}  =  0, &
\end{aligned}
\right. 
\end{equation}
in a domain $\Omega$ {\em with boundaries, endowed with a  no-slip boundary
condition.} Mathematically, the main difficulty is the lack of uniform bounds  on the
vorticity field, as the viscosity $\nu$ goes to zero. In terms of fluid
dynamics, this corresponds to a boundary layer phenomenon near $\pa
\Omega$.

\medskip
A natural approach to describe this boundary layer is to look for  a double-scale
asymptotics, with a parabolic scaling in the normal direction. Consider the case  $\Omega
\subset \R^2$. At least locally, any
point $\x$ in a neighborhood of $\pa \Omega$ has a unique decomposition  
$$\x \: = \:  y \,   \mathbf{n}(x) \:  + \:  \tilde{\x}(x), \quad \tilde{\x}\in
\pa \Omega,$$ 
where $y > 0$, $x$ is an arc length parametrization  of the boundary, and
$\mathbf{n}$ is the inward unit normal vector at $\pa \Omega$. The velocity
field can be written
$$   \u^\nu(t,\x) \: = \: u^\nu(t,x,y) \, {\bftau}(x) \: + \: v^\nu(t,x,y) \, \mathbf{n}(x),$$
where $(\bftau, \mathbf{n})$ is  the Fr\'enet frame.
It is then  natural to consider an approximation of the type:
\begin{equation} \label{ansatz}
\begin{aligned}
 u^\nu(t,x,y) \: & \approx   \: u^0(t,x,y) \: + \:  u^{BL}(t,x,y/\sqrt{\nu}), \\ 
v^\nu(t,x,y) \: & \approx   \: v^0(t,x,y) \: + \: \sqrt{\nu} \,  v^{BL}(t,x,y/\sqrt{\nu}),
\end{aligned}
\end{equation}    
where 
$$\u^0(t,x,y) \: = \:  u^0(t,x,y) \, {\bftau}(x) \: + \: v^0 (t,x,y) \,
\mathbf{n}(x)$$
satisfies the Euler equation with the no penetration condition, and 
 $$(u^{BL},v^{BL}) = (u^{BL},v^{BL})(t,x,Y)$$ 
describes a boundary layer corrector with  typical scale $\sqrt{\nu}$ in the normal
direction. It is slightly more convenient to introduce
$$ u(t,x,Y) \: := \: u^0(t,x,0) \: + \: u^{BL}(t,x,Y), \quad  v(t,x,Y) \:
 := \: Y \, \pa_y v^0(t,x,0) \: + \: v^{BL}(t,x,Y). $$
 Indeed, inserting this Ansatz  in the Navier-Stokes equations, we get formally 
\begin{equation} \label{Prandtl} 
\left\{
\begin{aligned}
 \pa_t u + u \pa_x u + v \pa_Y u  -  \pa^2_{Y} u    =  \left( \pa_t
u^0  + u^0 \pa_x u^0 \right)\vert_{y=0}, & \quad  Y > 0, \\ 
\pa_x u + \pa_Y v  = 0,   & \quad Y > 0, \\
(u,v)    = 0, & \quad Y = 0, \\
\lim_{Y \rightarrow +\infty}u = u^0\vert_{y=0}. &  
\end{aligned}   
\right.
\end{equation}
These are the so-called  Prandtl equations, derived by Ludwig Prandtl
\cite{Pra:1904}. Note that the curvature of the boundary does not appear explicitly
in the system. It is however involved in (\ref{Prandtl})  through the Euler field, and
through the interval of definition of 
the arc length parametrization $x$. Up to our knowledge, all studies deal
with one of the three following cases: $x \in \R$, $x \in \T$, or
$0 < x < L$, supplemented  with a condition on $u$ at $x=0$. The first
and second choices are convenient to describe  phenomena that are local in
$x$. The case  $x \in \T$ may also model  the outside of a bounded convex
obstacle. Finally, the third configuration is adapted to  the
spreading of a flow around a thin obstacle, where  $x=0$ corresponds to the tip
of the obstacle.

\medskip
Although this formal asymptotics is very natural, its validity is not clear.
As emphasized by physicists, including Prantdl himself, it may 
not hold uniformly in space and time. One reason is the
so-called {\em boundary layer separation}, which is observed for flows around
obstacles, see \cite{Guy:2001}. Nevertheless, the description \eqref{ansatz} fits
with  many experiments, upstream from the separation zone. In any
case, to understand the relevance and limitations of the Prandtl model is a
crucial issue.

\medskip
From the mathematical point of view, one must address two problems:   
\begin{enumerate}
\item The well-posedness of the Prandtl equation.
\item The justification of the expansion \eqref{ansatz}.   
\end{enumerate}
{\em These two problems depend crucially  on the choice of the underlying
functional spaces, especially on the regularity that is  required in the tangential
variable $x$.} Indeed, the main mathematical difficulty 
 is the  lack of control of the $x$ derivatives.
 For example, $v$ is recovered in  \eqref{Prandtl} through the divergence
 condition,  and in terms of $x$-regularity,  behaves broadly like $\pa_x
 u$. This loss of one derivative is not balanced by any horizontal
 diffusion term, so that standard energy estimates do not apply.

\medskip  
Within spaces of  functions that are analytic  in $x \in \R$, $Y \in \R^+$,  Sammartino and
Caflish  have overcome these problems, justifying locally in time 
 the boundary layer asymptotics \cite{Sam:1998a, Sam:1998b}. But for more ``realistic''
 functional settings,  the way  solutions of \eqref{NS} behave is still 
 poorly understood. Various instability mechanisms, that are  filtered
 out in an analytic framework, become  a huge source of trouble. For
 example, when the viscosity is small, the Navier-Stokes equation admits
 exponentially growing solutions which are both small-scale and isotropic
 in $x,y$. Their evolution is  lost in the anisotropic  Prandtl description.

\medskip
This remark was used by  Emmanuel Grenier in \cite{Gre:2000}, who relied on    
the so-called Rayleigh instability for inviscid flows to  show  that
{\em the asymptotics \eqref{ansatz} does not (always) hold  in the
  Sobolev space  $H^1$} (see \cite{Gre:2000} for a precise statement).
 However, the relevance of  this asymptotics   in $L^p$, or its relevance in the
 absence of Rayleigh instabilities, are still open issues.

\medskip
{\em Above all, the local in time  well-posedness of the  Prandtl equation
  for smooth (say Sobolev) initial
data has been so far an open problem}. Up to our knowledge, the Cauchy
problem has only been solved in two settings: 
\begin{description}
\item[i)] $x \in \R$,  with data  that are analytic in $x$ : see
  \cite{Sam:1998a,Lom:2003}  for more.
\item[ii)] $0 < x < L$, with data that are monotonic in $y$ : see
  \cite{Ole:1999,Xin:2004} for more. 
\end{description}
One may  also cite article \cite{E:1997}, in which blow up in time of some
smooth solutions is exhibited. Finally, let us mention the interesting work
\cite{Hon:2003}, in which the inviscid version of \eqref{Prandtl} is analyzed (no
$\pa^2_Y u$ in the equation). Interestingly, for a smooth initial data, 
 this  equation turns out to have an explicit solution  through the
 method of characteristics. In particular, starting from a smooth data, one
 recovers locally in time a smooth data. More precisely, there is only a
 finite loss of $x$-derivatives, so that the  Cauchy problem is
 (weakly) well-posed. We refer to \cite{Hon:2003} for all details. See
 also papers \cite{Gre:1999,Bre:1999} on the hydrostatic equations, that share some 
features with Prandtl equations. For more on Prandtl equations, see the
review \cite{E:2000}.

\medskip
On the basis of the inviscid  result, it seems reasonable to bet for
 well-posedness of the Prandtl equation \eqref{Prandtl} in Sobolev type
spaces. {\em The aim of this paper is to show that it is
actually linearly ill-posed in this framework}. As we shall see later on, the reason for
ill-posedness is a strong destabilization mechanism due to two ingredients:
viscosity, and critical points  in the base velocity profile. In particular, it does
not contradict the positive results obtained in the inviscid case and for
monotonic data. 

\medskip
Let us now describe  our results. We restrict ourselves to  $(x,Y) \in \T
\times \R^+$, and $\u^0 = 0$. To lighten notations, we write $y$ instead of
$Y$. The Prandtl equation comes down  to 
\begin{equation} \label{Prandtl2} 
\left\{
\begin{aligned}
 \pa_t u + u \pa_x u + v \pa_y u  -  \pa^2_{y} u    =  0,  & \quad  \mbox{ in } \: \T
\times \R^+. \\ 
\pa_x u + \pa_y v  = 0,   & \quad \mbox{ in } \: \T
\times \R^+, \\
(u,v)\vert_{y=0}   = (0,0),  & \quad \lim_{y\rightarrow +\infty} u = 0.
\end{aligned}   
\right.
\end{equation} 
Let $u_s = u_s(t,y)$ a smooth solution of the heat equation 
\begin{equation} \label{heat}
 \pa_t u_s - \pa^2_y u_s = 0, \quad u_s\vert_{y=0}=0, \quad
u_s\vert_{t=0} = U_s,
\end{equation}
with good decay as $y \rightarrow +\infty$.   Clearly, the shear velocity profile  $(u_s,v_s) = (u_s(t,y), 0)$ satisfies 
the  system  \eqref{Prandtl2}. We consider the linearization around
$(u_s,v_s)$, that is 
\begin{equation} \label{linear} 
\left\{
\begin{aligned}
 \pa_t u + u_s \pa_x u + v \pa_y u_s  -  \pa^2_{y} u    =  0,  & \quad  \mbox{ in } \: \T
\times \R^+. \\ 
\pa_x u + \pa_y  v  = 0,   & \quad \mbox{ in } \: \T
\times \R^+, \\
(u,v)\vert_{y=0}   = (0,0), & \quad  \lim_{y\rightarrow+\infty} u = 0.
\end{aligned}   
\right.
\end{equation} 
{\em We wish to study well-posedness properties of  \eqref{linear}, for a certain
class  of velocities  $u_s$.} In this view, we introduce the following
  functional spaces:
$$W^{s,\infty}_\alpha(\R^+) \: := \:  \left\{ f = f(y), \quad e^{\alpha y}f
\in W^{s,\infty}(\R^+) \right\}, \quad \forall \alpha, s \ge 0, $$
with  $\: \| f \|_{W^{s,\infty}_\alpha} \: := \: \| e^{\alpha y} f
\|_{W^{s,\infty}}, \: $ and  
$$E_{\alpha,\beta} \: := \:  \left\{ u = u(x,y) = \sum_{k \in \Z}
\hat{u}^k(y) e^{i k x}, \quad
 \|  \hat{u}^k
\|_{W^{0,\infty}_\alpha} \, \le \, C_{\alpha,\beta} \ e^{-\beta |k|}, \: \forall k
\right\}, \quad \forall \alpha, \beta >  0, $$ 
with $\: \| u \|_{E_{\alpha,\beta}} \: := \: \sup_{k} \, e^{\beta |k|} \,  \|  \hat{u}^k
\|_{W^{0,\infty}_\alpha}.$

\medskip
Note that   the functions of  $E_{\alpha,\beta}$ have analytic
regularity in $x$. They have only  $L^\infty$ regularity in $y$, with an
exponential weight. More regularity in $y$ could be
 considered as well. Let $\alpha, \beta > 0$.  We prove in the appendix the following result: 
\begin{proposition} \label{propanalytic} {\bf (Well-posedness in the
    analytic setting)} 

\smallskip
Let $u_s \in C^0\bigl(\R_+; \,W^{1,\infty}_{\alpha}(\R_+)\bigr)$.  
There exists $\rho > 0$ such that: for all $T$
with $\beta - \rho T > 0$, and all $u_0 \in E_{\alpha,\beta}$, the linear
equation \eqref{linear} has a unique solution 
$$\displaystyle u \in C\left([0,T); \,  E_{\alpha, \beta-\rho T}\right),
  \quad u(t,\cdot) \in  E_{\alpha, \beta-\rho t}, \quad   
u\vert_{t=0} = u_0.$$ 
\end{proposition}     
In short, the Cauchy problem for \eqref{linear} is locally well-posed in the
analytic setting.  We shall  denote  
$$ T(t,s) u_0 \: := \: u(t, \cdot) $$
where $u$ is the solution of \eqref{linear} with $u\vert_{t=s} = u_0$.  
As the spaces $E_{\alpha,\beta}$ are dense in the  spaces
$$H^m \: := \: H^m(\T_x,W^{0,\infty}_\alpha(\R^+_y)), \quad m \ge 0,$$
this makes sense to introduce the following notation: for all $T \in {\cal
  L}(E_{\alpha,\beta}, E_{\alpha,\beta'})$,
$$ \| T \|_{{\cal L}(H^{m_1},H^{m_2})}  \: = \: \sup_{u_0 \in
  E_{\alpha,\beta}} \, \frac{\|T  u_0\|_{H^{m_2}}}{\| u_0
  \|_{H^{m_1}}},$$
that belongs to $\R_+ \cup \{+\infty\}$. In particular, it is infinite when
$T$ does not extend to a bounded operator from  $H^{m_1}$ to $H^{m_2}$. The main  result of our  paper is
\begin{theorem} \label{theorem} {\bf (Ill-posedness in the Sobolev setting)}
\begin{description}
\item[i)] Let $\displaystyle u_s \in C^0\bigl(\R_+; \,W^{4,\infty}_{\alpha}(\R_+)\bigr)
  \cap C^1\bigl(\R_+; \,W^{2,\infty}_{\alpha}(\R_+)\bigr)$. 
\smallskip
Assume that the initial velocity  has a non-degenerate critical point over
$\R^+$. Then, there exists $\sigma > 0$, such that for all $\delta > 0$, 
\begin{equation*}
 \hspace{1.5cm} \sup_{0 \le s \le t \le \delta} \: \|e^{-\sigma (t-s)\sqrt{|\pa_x|}} \,
 T(t,s) \|_{{\cal L}(H^{m},H^{m-\mu})}  \: = 
\: + \infty, \quad \forall m \ge 0, \: \mu \in [0,1/2). 
\end{equation*}
\item[ii)] Moreover, one can find  solutions $u_s$ of \eqref{heat} and
  $\sigma > 0$ such that:  for all $\delta > 0$,
 \begin{equation*}
\sup_{0 \le s \le t \le \delta} \: \|e^{-\sigma (t-s)\sqrt{|\pa_x|}} \,  T(t,s) \|_{{\cal L}(H^{m_1},H^{m_2})}  \: =
  \: + \infty, \quad \forall m_1,m_2 \ge 0.  
\end{equation*}
\end{description}
\end{theorem}
This theorem expresses strong linear ill-posedness of the Prandtl equation
in the Sobolev framework. It is a  consequence of an instability process, that holds at high tangential
frequencies. We will show that  some perturbations with tangential
frequency $k \gg 1$ grow like $e^{\sqrt{k}t}$.

\medskip 
The outline of the paper is as follows. Section \ref{sec1} gives a formal
description of the instability mechanism. It relies
on an asymptotic analysis of \eqref{linear}, in the
high tangential frequency limit. Thanks to this analysis, we show that
 ill-posedness for the  PDE \eqref{linear} comes down to a ``spectral
 condition''  for  a reduced ODE, namely: 

\medskip
\noindent
{\em (SC) There exists $\tau \in \C$ with $\Im \tau < 0$, and a solution $W = W(z)$  of
 
\begin{equation}  \label{ODE}
(\tau - z^2)^2 \frac{d}{dz} W \: + i \:
\frac{d^3}{dz^3} \left( (\tau - z^2) W \right) = 0, 
\end{equation}
such that $\displaystyle \lim_{z \rightarrow -\infty} W = 0, 
\lim_{z \rightarrow +\infty} W = 1$.}

\medskip  
\noindent
This spectral condition is studied in section \ref{sec2}, and shown to be
satisfied. On these grounds, we  prove theorem \ref{theorem}, {\it cf}
section \ref{sec3}. We end up the paper with  numerical computations, which
emphasize that our  instability mechanism is effective.

\section{The instability mechanism} \label{sec1}
In this section, we  describe the destabilization of system
\eqref{linear}, leading to the ill-posedness theorem. As we shall see, it
takes place at high tangential frequencies, say $O(1/\eps)$, and has a typical  time 
$O(\sqrt{\eps})$. At this timescale, the time dependence of the base velocity
$(u_s(t,y),0)$ will not play an important role. Thus, to understand the
instability mechanism, we can consider the simpler equation
\begin{equation} \label{linear2} 
\left\{
\begin{aligned}
 \pa_t u + U_s \pa_x u + v  U_s'  -  \pa^2_{y} u    =  0,  & \quad  \mbox{ in } \: \T
\times \R^+. \\ 
\pa_x u + \pa_y  v  = 0,   & \quad \mbox{ in } \: \T
\times \R^+, \\
(u,v)\vert_{y=0}   = (0,0). &
\end{aligned}   
\right.
\end{equation} 
 Handling of the real equation, that is with $u_s$ instead of $U_s$, will
require minor modifications, to be achieved in section \ref{sec3}.

\medskip
System \eqref{linear} has constant coefficients in $t$ and $y$, so that we
can perform a Fourier analysis: we look for solutions in the form
\begin{equation} \label{Fourier}
 u(t,y) \: = \:  e^{i k (\omega(k) t + x)} \hat{u}^k(y), \quad v\: = \:  k \,
e^{i k (\omega(k) t + x)} \hat{v}^k(y), \quad k > 0.  
\end{equation} 
As we are interested in high tangential frequencies, we denote $\eps \:
:= \: 1/k \ll 1$, and write $\omega(\eps)$ instead of $\omega(k)$, $\: u_\eps(y),
v_\eps(y)$ instead of $\hat{u}^k(y),
\hat{v}^k(y)$.  The divergence condition yields $v_\eps'(y) \: = \: -i
u_\eps(y)$. Using this relation in the first equation in \eqref{linear2}, one ends
up with 
\begin{equation} \label{linear3} 
\left\{
\begin{aligned}
\left( \omega(\eps) +  U_s \right) v_\eps' \: - \:    U_s'  v_\eps  \: + \:
i \, \eps v_\eps^{(3)}  \:  = \:   0,  & \quad y > 0, \\ 
v_\eps\vert_{y=0}  \: = \: v_\eps'\vert_{y=0} = 0. &
\end{aligned}   
\right.
\end{equation} 
Thus, the high frequency limit $\eps \rightarrow 0$ in variable $x$  yields a singular
perturbation problem in variable $y$. To investigate this problem, one must
first consider the inviscid case $\eps=0$.

\subsection{The inviscid case}
When $\eps =0$, one can  {\it a priori} only retain the impermeability
condition. The appropriate problem  is 
\begin{equation} \label{inviscidlinear}
\left( \omega +  U_s \right) v' \: - \:    U_s'  v   \:  = \:   0, 
\quad  y > 0, \quad  v\vert_{y=0} = 0.
\end{equation}
This spectral problem, as well as the corresponding evolution equation,
have been studied exhaustively in \cite{Hon:2003}. Clearly, there are  non-trivial solutions
if and only if  $\omega$ belongs to the range of $-U_s$. 
Moreover, the couples
$$ \omega_a \: = \: -U_s(a), \quad v_a \: = \:   H(y-a) \, \left( U_s - U_s(a) \right),
\quad a >  0$$ 
where $H$ is the Heaviside function, satisfy \eqref{inviscidlinear}.
Note that the regularity of $v_a$ depends on the choice of $a$. When $a$ is
a critical point, it belongs to  $W_\alpha^{2,\infty}(\R^+)$  with a
discontinuous second derivative. Otherwise, it is only in $W^{1,\infty}_\alpha(\R^+)$, with a
 discontinuous first derivative. Luckily enough, the additional boundary
 condition $v_a'\vert_{y=0} = 0$ is also satisfied. 

\subsection{The viscous perturbation}
When $\eps$ is not  $0$, the inviscid eigenelements $\omega_a, v_a$ 
do not solve \eqref{linear3}. All boundary conditions
are satisfied, {\it cf} the above remark, but the equation is not.  First, there is a
$O(\eps)$ remaining term for $y >  a$. More importantly, $v_a$ is not
smooth at $y=a$, whereas a solution of this parabolic equation
should be. 

\medskip
Nevertheless, at least  if   $a$ is a non-degenerate critical point,
there is  an approximate solution  near $(\omega_a,v_a)$. 
We shall establish this  rigorously in section  \ref{sec3}.  We just give here a 
formal expansion. It reads
\begin{equation} \label{ansatz2}
\left\{
\begin{aligned}
\omega(\eps) \: & \sim \: \omega_a \: + \: \eps^{1/2} \tau, \\
v_\eps(y) \: & \sim  \: \, v_a \, + \,
\eps^{1/2} \tau \,  H(y-a)  \: + \:  \eps^{1/2} \ V\left( \frac{y-a}{\eps^{1/4}} \right),
\end{aligned}
\right.
\end{equation} 
where $\tau \in \C$, and $V = V(z)$ quickly tends to zero as $z \rightarrow
\pm \infty$. Note that the approximation of $v_\eps$ has  two
parts: the  ``regular'' part  
$$v_\eps^{reg}(y) \: = \:  H(y-a)  \, \left( U_s(y) \, - \, U_s(a) \, + \, \eps^{1/2} \tau
\right)$$
and the ``shear layer part''
$$ v_\eps^{sl}(y) \: = \:  \eps^{1/2} \ V\left( \frac{y-a}{\eps^{1/4}}
\right). $$
For $\omega(\eps) = -U_s(a) \: + \: \eps^{1/2} \tau$, the function  $v_\eps^{reg}$ solves
\eqref{linear3} up to $O(\eps)$,  away from the critical point
$y=a$. However, it has a jump at $y=a$, together with its second
derivative.
 The role of the shear layer
$v_\eps^{sl}$, which concentrates near $y=a$, is to cancel these discontinuities. Still formally, we
obtain the system satisfied by the profile  $V$:   
\begin{equation} \label{SL}
\left\{
\begin{aligned}
& \left(\tau + U_s''(a) \frac{z^2}{2}\right) V' \: - \:  U_s''(a)  \, z \, V \: +
\: i  \,  V^{(3)} \: = \: 0, \quad  \: z \neq 0, \\ 
& \left[ V \right]_{\vert_{z=0}} \: = \:  - \tau, \quad  \left[ V' \right]_{\vert_{z=0}} \: =
\: 0, \quad  \left[ V'' \right]_{\vert_{z=0}} = - U''(a), \\
& \lim_{\pm\infty} V \: = \: 0. 
\end{aligned}
\right.
\end{equation}
Let us point out that this system is {\it a priori} overdetermined, as jump
and boundary conditions provide too many constraints. This justifies the introduction of the
parameter $\tau$ in the Ansatz \eqref{ansatz2}. As we shall see below,
there is  a $\tau$ for  which system
\eqref{SL} has a solution. {\em Moreover, $ \Im \tau$ is negative.
 Hence, back to the Fourier representation \eqref{Fourier}, the
  $k$-th mode will grow in time like $e^{\sqrt{k}t}$. This is the key of
the instability mechanism.}  

\medskip
To see how the condition (SC) of the introduction steps in, we need a few
rewritings. First, $\tau + U_s''(a) \frac{z^2}{2}$ satisfies
the equation in \eqref{SL}. We therefore introduce 
$$ \tilde V(z) \: = \: V(z) \: + \: {\bf 1}_{\R_+} \left(  \tau + U_s''(a)
\frac{z^2}{2} \right) $$
which leads to
\begin{equation*} 
\left\{
\begin{aligned}
& \left(\tau + U_s''(a) \frac{z^2}{2}\right)  \tilde V' \: - \:  U_s''(a)
  \, z \,  \tilde V \: + 
\: i  \,  \tilde V^{(3)} \: = \: 0, \quad z \in \R, \\ 
& \lim_{-\infty}  \tilde V \: = \: 0, \quad 
  \tilde V \: \sim_{+\infty} \: \tau + U_s''(a) \frac{z^2}{2}. 
\end{aligned}
\right.
\end{equation*}
Then, we introduce $W$  such that 
$$\tilde V \: = \:  \left( \tau \, +
\, U_s''(a) \, z^2/2 \right) \,  W.$$
 We get:
\begin{equation*} 
\left\{
\begin{aligned}
& \left(\tau + U_s''(a) \, z^2/2 \right)^2   \frac{d}{dz} W \: +
\: i \, \frac{d^3}{dz^3} \left(  \left( \tau + U_s''(a)
\, z^2/2 \right) \,  W \right)  \: = \: 0, \\ 
& \lim_{-\infty}   W \: = \: 0, \quad \lim_{+\infty}
  W \: = \: 1.
\end{aligned}
\right.
\end{equation*}
Finally, we perform the change of variables
$$ \tau \: = \: \frac{1}{\sqrt{2}}|U''(a)|^{1/2} \, \tilde \tau, \quad z \: =
  2^{1/4} \, |U''(a)|^{-1/4}  \tilde z. $$
Dropping the tildes leaves us with the reduced ODE 
\begin{equation*} 
\left\{
\begin{aligned}
& \left(\tau + \mbox{ sign}(U_s''(a)) \, z^2 \right)^2   \frac{d}{dz} W \: +
\: i\, \frac{d^3}{dz^3} \left(  \left( \tau + \mbox{ sign}(U_s''(a))
\, z^2 \right) \,  W \right)  \: = \: 0, \\ 
& \lim_{-\infty}   W \: = \: 0, \quad \lim_{+\infty}
  W \: = \: 1.
\end{aligned}
\right.
\end{equation*}
If $U''_s(a) < 0$, it is exactly the system in (SC). If on the contrary $U''_s(a) > 0$, and if $(\tau, W)$ satisfies
the system in  (SC), then $\displaystyle (\tau := -\overline{\tau}, \, W := \overline{W})$ satisfies the
above system. In both cases, back to the original system \eqref{SL}, condition
(SC) gives a solution $(\tau, V)$ with $\Im \tau < 0$.   
In particular,  this $\sqrt{\eps}$ correction to the eigenvalue is
a source of strong instability, leading to  ill-posedness.

\medskip
The proof of Theorem \ref{theorem}, which is based on this formal
shear layer phenomenon, is postponed to section \ref{sec3}. In the next paragraph, we
focus on condition (SC), and prove that it is satisfied.

\section{The spectral condition (SC)} \label{sec2}
We need to study the existence of  heteroclinic orbits for the ODE
\eqref{ODE}. Note that $W=1$ is  a solution.  Equation \eqref{ODE} can be written
as a second order equation in $X=W'$:
\begin{equation} \label{ODE2}
 i(\tau - z^2) X'' \: - \: 6  i \,  z \,  X' \: + \: \left((\tau - z^2)^2 - 6i\right) X \: =
\: 0. 
\end{equation}
To show that (SC) holds, we proceed in three steps. 

\bigskip
{\em Step 1.}
We consider an auxiliary eigenvalue problem:
\begin{equation} \label{eigenpb}
A u \: := \: \frac{1}{z^2+1} u'' \: + \: \frac{6z}{(z^2+1)^2} u' \: + \:
\frac{6}{(z^2+1)^2} u \: = \: \alpha \, u.
\end{equation}
For its study, we introduce the weighted spaces 
\begin{align*}
{\cal L}^2 \: & := \: \left\{ u \in L^2_{loc}, \:  \int_\R (z^2 + 1)^4 | u |^2  < +\infty
\right\}, \\
 {\cal H}^1 \: &  := \: \left\{ u \in H^1_{loc}  \:  \int_\R (z^2 + 1)^4 | u
 |^2  \: + \: \int_\R (z^2 + 1)^3 | u' |^2 < +\infty \right\} \:
 \hookrightarrow {\cal L}^2. 
\end{align*}
with their obvious Hilbert norms. We see  $A$ as an operator from $D(A) \:
:= \: \left\{ u \in {\cal H}^1, \; Au \in {\cal L}^2 \right\}$ into ${\cal
  L}^2$. {\em Our goal is to show that $A$ has a positive eigenvalue.}

\medskip
By standard arguments, the domain $D(A)$ is dense in ${\cal L}^2$. Moreover, for any
 $u$ in $D(A)$, there is a sequence $u_n$ of smooth functions with compact
 support, such that $u_n \rightarrow u$ in ${\cal H}^1$ and $Au_n \rightarrow Au$ in ${\cal
   L}^2$. Integration by parts and use of this density property give
 easily  that  $A$ is symmetric, {\it i.e.} 
\begin{equation} \label{A1}
\forall u,v \in D(A), \quad (A u \, | \,  v)_{{\cal L}^2} \: = \:  (A v \,
| \,  u)_{{\cal L}^2} 
\end{equation}
and that for $\lambda$ large enough
\begin{equation} \label{A2}
((\lambda - A) u \, | \,  u)_{{\cal L}^2} \: =  \:    \lambda \int_\R (z^2 + 1)^4 | u
 |^2   - 6 \int_\R (z^2 + 1)^2 | u |^2 \: + \:  \int_\R (z^2 + 1)^3 | u' |^2
 \: \ge \: \frac{1}{2} \| u \|^2_{{\cal H}^1}. 
\end{equation}
Then, the coercivity condition \eqref{A2} allows to apply  the Lax-Milgram
lemma. It  implies the invertibility of $\lambda -A$, with 
$$ \|\left( \lambda-A \right)^{-1} f \|_{{\cal L}^2} \: \le \: \|\left(
\lambda-A \right)^{-1} f \|_{{\cal H}^1} \: \le \: C \, \| f \|_{{\cal
    L}^2}. $$
Moreover, from \eqref{A1}, $\left( \lambda-A \right)^{-1}$ is selfadjoint,
and so is $A$. 

\medskip
First, let us prove  that $A$ has positive spectrum. To do so, we claim
that it is enough to find $u \in  D(A)$ with $(A u \, | \, u)_{{\cal L}^2}
> 0$. Indeed,   suppose {\it a contrario} that $\sigma(A)$ is contained in
$\R_-$. Then, by the spectral theorem, 
$$ \forall \alpha > 0, \quad  \| (A - \alpha)^{-1} \| =
\frac{1}{d(\alpha,\sigma(A))} \: \le  \: \alpha^{-1}. $$  
We  deduce:  for all $u \in D(A)$ 
$$ \| u \|^2_{{\cal L}^2} \: \le \:  \alpha^{-2} \,  \|(A - \alpha) u \|^2_{{\cal L}^2} $$ 
Expanding the scalar products, we obtain
$$  0 \: \le \:  \alpha^{-2}  \,  \|A u \|^2_{{\cal L}^2} \: - \: 2
\alpha^{-1} (A u \, | \, u)_{{\cal L}^2} $$
In the limit $\alpha \rightarrow + \infty$, we get  $(A u \, | \, u)_{{\cal L}^2} \le
0$ for all $u \in D(A)$.  This proves our claim.
From there, we simply  take $u = e^{-2 z^2}$. A straightforward computation gives
$$ (A u \, | \,  u ) \: = \: \frac{439}{512} \, \sqrt{\pi} \: > \: 0, $$
and so $\sigma(A)$ has a  positive subset.

\medskip
It remains to exhibit a positive eigenvalue inside this positive subset of
the spectrum. We remark that the operator $A$ can be split into 
$$ A \: = \:  A_1 \: + \: A_2, \quad A_1 u \: := \: \frac{1}{z^2+1} u'' \: +
\: \frac{6z}{(z^2+1)^2} u', \quad A_2 u \: :=  \: \frac{6}{(z^2+1)^2} u. $$
On one hand,  the operator $A_1$  is negative, and by Lax-Milgram Lemma, $A_1 -
\lambda$ is invertible for any $\lambda > 0$. Thus,  $\sigma(A_1) \subset
\R_-$. On the other hand, let $u_n$ a sequence with $u_n$ and $A_1 u_n$
bounded in ${\cal L}^2$. This implies that $u_n$ is bounded in ${\cal
  H}^1$, and so has a convergent subsequence in
$L^2_{loc}$. Moreover, $|u_n|^2$ is equi-integrable over $\R$. Finally, it
implies that  $A_2 u_n$ has a convergent subsequence in ${\cal L}^2$, which means that  
 $A_2$ is  $A_1-$compact. Hence, the essential spectra of $A$ and $A_1$
are the same, see \cite{Kat:1995}. In particular, the positive part of $\sigma(A)$
is made of isolated eigenvalues with finite multiplicity.         
Eventually, we state:  {\em there exists $\alpha >
  0$, and $u$ in $D(A)$ satisfying \eqref{eigenpb}.}

\bigskip
{\em Step 2.} We wish to convert  the eigenelements $(\alpha,u)$ of the
previous step  into an appropriate solution $(\tau,X)$ of \eqref{ODE2}.
 We set $\tilde \tau = -\alpha^{1/2}$, and $\tilde z = \alpha^{-1/4}z$, $\:
 Y(\tilde z) = u(z)$. Dropping the tildes, we obtain  a solution of 
\begin{equation} \label{ODE3}
 (\tau - z^2) Y'' \: - \: 6   \,  z \,  Y' \: + \:
  \left((\tau - z^2)^2 - 6\right) Y \: =
\: 0. 
\end{equation}
By a classical bootstrap argument,  $Y$ is smooth. Moreover, it inherits
from $u$ its integrability properties at infinity. Actually, the
behaviour of $Y$ can be further specified, as shown in: 
\begin{proposition} \label{behavODE}
 The function $Y$ admits a unique extension, still denoted by $Y$, that
 is holomorphic in $z$ and satisfies \eqref{ODE3}
  in the simply connected  domain 
$$U_\tau \: := \:\C\setminus \left( \left[-i \infty, \,
    -i|\tau|^{1/2}\right] \, \cup \,  \left[i |\tau|^{1/2}, \, +i\infty\right]\right).$$
 Moreover, in the sectors  $\arg  z  \in (-\pi/4+\delta, \pi/4-\delta)$ and  $ \arg
  z \in (3\pi/4+\delta, 5\pi/4-\delta)$, $\delta > 0$, it satisfies the inequality
$$ | Y(z) | \: \le \: C_\delta \, \exp(-z^2/4). $$
\end{proposition}

\medskip
{\em Proof.} This proposition  follows from the general theory of
ODE's with holomorphic coefficients. The existence of a
holomorphic solution is well-known, because the coefficient $
\tau - z^2$ does not vanish on $U_\tau$. As regards the inequality,
 we rewrite  equation \eqref{ODE2}  as the first order system:  
\begin{equation} \label{firstorder}
\frac{d}{dz} {\cal Y} \: = \:  z {\cal A}( z) {\cal Y}, \quad {\cal
  Y} \: = \: \left( \begin{smallmatrix} Y  \\ z^{-1}
  \frac{d}{dz} Y
 \end{smallmatrix} \right), 
  \quad {\cal A}(z) \: = \: \left( \begin{smallmatrix} 0 & 1 \\ \frac{6 -
      (\tau - z^2)^2}{ z^2 (\tau - z^2)}  & \frac{6}{\tau - z^2} -
    \frac{1}{z^2} \end{smallmatrix} \right).
\end{equation}
In particular,  ${\cal A}$ is holomorphic at infinity, with
$\displaystyle  {\cal A}(\infty) \: = \: \left( \begin{smallmatrix} 0 & 1 \\ 1  &
  0 \end{smallmatrix} \right)$.
It has two distinct eigenvalues $\pm 1$, with eigenvectors 
$\left( \begin{smallmatrix} 1  \\ \pm 1   \end{smallmatrix} \right).$

\medskip
Hence, we can apply  \cite[Theorem 5.1 p163]{Cod:1955}:  
in any closed sector $S$ inside which $\Re
z^2$ does not cancel, there exists solutions  ${\cal Y}_\pm$ (depending {\it a priori}
on S)  with the following asymptotic behaviour as $|z| \rightarrow  +\infty$:
\begin{equation} \label{ansatz3}
  {\cal Y}_{\pm}   \: \sim  \:  \left(  \sum_{i\ge 0}   {\cal
    Y}^i_\pm \, z^{\alpha_{\pm}-i}  \right) \, 
e^{P_{\pm}(z)}, \:  
\end{equation}
where $\alpha_{\pm}$ is a complex constant, and  
$P_{\pm}(z)$ is a polynomial of degree 2. Moreover, the  leading term of
$P_{\pm}$ is $\pm \frac{z^2}{2}$. 
 Following the same scheme of proof, we get in the present case:
\begin{equation} \label{ansatz4}
{\cal Y}^0_\pm \: = \:  \left( \begin{smallmatrix} 1
  \\ \pm 1  \end{smallmatrix} \right), \quad  \alpha_\pm \: = \: 
  -\frac{1}{2} \left( \pm \tau    \, +\,  7 \right), \quad P_\pm \: = \: \pm z^2/2.
\end{equation}     
As our solution $Y$ is integrable over  $\R$, it is necessarily
proportional to the decaying solution.  The bounds in proposition \ref{behavODE} follow. 

\bigskip
Now, as $Y$ is defined on $U_\tau$, we can perform the complex change
of variable:
$$ \tilde{z} \: := \: e^{i\pi/8} z, \quad \tilde \tau \: := \:  e^{i\pi/4} 
\tau $$ 
Note that $\tilde \tau$ has a negative imaginary part. Moreover, for
$\tilde z$ real, the original variable $z$ belongs to the sectors  
$\displaystyle \arg  z  \in (-\pi/4+\delta, \pi/4-\delta)$ or
$\displaystyle  \arg
  z \in (3\pi/4+\delta, 5\pi/4-\delta)$, with $\delta =
  1/16$. By the proposition \ref{behavODE}, the function $ X(\tilde z) \: = \: Y(z)$ satisfies the
  estimate $\displaystyle |X(\tilde z)| \: \le \: C \, \exp(-1/4 \, \tilde z^2)$. Finally, dropping
  the tildes yields a solution $\tau, X$ of \eqref{ODE2}, where $X$ decays
  at infinity. This concludes step 1. 

\bigskip
{\em Step 3.} To deduce from the previous step that  (SC) holds, it is
enough that  $\int_\R X(z) dz$ be non zero. If so, one can  consider  
$$W(z) \: := \:  \left( \int_\R X(z') dz' \right)^{-1} \int_{-\infty}^z X(z')
  \, dz',$$ 
which  clearly satisfies all requirements.

\medskip
Let us assume {\it a contrario} that $\int_\R X(z) dz =0$. Then, the
function
$$ V(z) \: := \:  (\tau - z^2) \int_{-\infty}^z X(z') \, dz' $$
is a solution of 
\begin{equation*} 
\left(\tau - z^2 \right)   V' \: + \: 2 \, z \,  V \: +
\: i  \,  V^{(3)} \: = \: 0, 
\end{equation*}
which decays exponentially as $z$ goes to  $\pm \infty$, together with all its
derivatives. Differentiation of the equation gives
\begin{equation*} 
\left(\tau - z^2 \right)   V'' \: + \: 2 \,   V \: +
\: i  \,  V^{(4)} \: = \: 0, 
\end{equation*}
Then, we multiply by $\overline{V''}$, that is the complex conjugate of
$V''$, and integrate over $\R$. Simple integrations by parts yield:
$$ \int_\R (\tau - z^2)  | V''|^2 \: - \: 2 \int_\R | V' |^2 - i \int_\R |
V^{(3)}|^2 \: = \: 0.$$   
The imaginary part of this identity yields  
$$ \Im \tau  \int_\R  | V''|^2  = \int_R | V^{(3)}|^2 $$
which contradicts the fact that $\Im \tau < 0$. {\em Thus, the condition  (SC) is
satisfied.} 

\section{Proof of ill-posedness} \label{sec3}
Theorem \ref{theorem} will be deduced from the formal analysis
of section \ref{sec1}. This analysis was performed on
\eqref{linear2}, in which possible time variations of $u_s$were  neglected.
 To account for the original system \eqref{linear} will
require a few modifications, notably in the choice of the approximation
\eqref{ansatz2}.  We will distinguish between the parts i) and ii) of the theorem.

\subsection{Ill-posedness for general  $u_s$}
Let $u_s$ satisfying the assumptions of  part i).   Let $a$ be the non-degenerate
critical point of $u_s\vert_{t=0} = U_s$. For the sake of brevity, we
consider the case $U"_s(a) < 0$, the other one being strictly similar. 
The  differential equation
$$ \pa_t\pa_y u_s(t,a(t)) \: + \: \pa_y^2 u_s(t,a(t)) \,a'(t) \: = \: 0, \quad
a(0) = a $$
defines for small time $t < t_0$ a non-degenerate  critical point $a(t)$ of
$u_s(t,\cdot)$. Let then  $\tau,W$ given by condition (SC). We set
$$V \: := \:  (\tau - z^2) W \: - \: {\bf 1}_{\R_+} \left( \tau - z^2 \right). $$  
In the light of section \ref{sec1}, we introduce, for $\eps > 0$ and $t
< t_0$: 
$$ \omega(\eps,t) \: := \:  - u_s(t,a(t)) \: + \: 
\frac{\eps^{1/2}}{\sqrt{2}} | \pa^2_yu_s(t,a(t))|^{1/2} \tau $$
 as well as  the ``regular'' velocity
$$v_\eps^{reg}(t,y) \: := \:  H(y-a(t)) \left( u_s(t,y) - u_s(t,a(t)) \:
  +\:   
\frac{\eps^{1/2}}{\sqrt{2}} | \pa^2_y u_s(t,a(t))|^{1/2} \tau \right), $$  
and the shear layer velocity
$$ v_\eps^{sl}(t,y) \: := \: \frac{\eps^{1/2}}{\sqrt{2}} \varphi(y-a(t)) \,
 | \pa^2_yu_s(t,a(t))|^{1/2}  \, 
V\left(|\pa^2_y u_s(t,a(t))|^{1/4}\frac{(y-a(t))}{(2 \eps)^{1/4}}\right), $$
where $\varphi$ is a smooth truncation function near $0$. 
We then consider the following  velocity field:
\begin{align*}
u_\eps(t,x,y) \:  & := e^{i \eps^{-1} x}  U_\eps(t,y), \quad U_\eps(t,y) \:
=\: i \, e^{i \eps^{-1} \int_0^t \omega(\eps,s) ds}   \, \pa_y \left(v^{reg}_{\eps}(t,y)  \: +
  \: v^{sl}_{\eps}(t,y)\right), \\
 v_\eps(t,x,y) \: & := \: e^{i \eps^{-1} x}  V_\eps(t,y), \quad V_\eps(t,y)
 \: = \: \eps^{-1} \, e^{i \eps^{-1}  \int_0^t \omega(\eps,s) ds}  \, \left(v^{reg}_{\eps}(t,y) \: +
  \: v^{sl}_{\eps}(t,y)\right).
\end{align*}
In order to have a field that is $2\pi$-periodic in $x$ and growing in time, we
 take  $\eps := \frac{1}{n}, \: $ with $\: n \in
  \N_*$.
One verifies easily that  $u_\eps = e^{i \eps^{-1} x} \, U_\eps(t,y)$  is
analytic in $x$, and $W^{2,\infty}$ in $t,y$. Moreover, we have the bounds 
\begin{equation} \label{expo1}
 c \, e^{\frac{\sigma_0 t}{\sqrt{\eps}}} \: \le \:   \| U_\eps(t,\cdot)
\|_{W^{2,\infty}_\alpha} \: \le \: C  \,e^{\frac{\sigma_0 t}{\sqrt{\eps}}}, 
\end{equation}
for positive constants $c, C$ and $\sigma_0$ that do not depend on $\eps$.

\medskip
Inserting the expression for $u_\eps, v_\eps$ into  the linearized Prandtl
equation \eqref{linear}, we obtain 
\begin{equation} \label{linear4} 
\left\{
\begin{aligned}
 \pa_t u_\eps + u_s \pa_x u_\eps + v_\eps \pa_y u_s  -  \pa^2_{y} u_\eps
 =  r_\eps,  & \quad  \mbox{ in } \: \T 
\times \R^+. \\ 
\pa_x u_\eps + \pa_y  v_\eps  = 0,   & \quad \mbox{ in } \: \T
\times \R^+, \\
(u,v)\vert_{y=0}   = (0,0). &
\end{aligned}   
\right.
\end{equation} 

The remainder term $r_\eps$ reads $r_\eps \: = \: e^{i \eps^{-1} x}
R_\eps(t,y)$, with   
\begin{align*}
 R_\eps(t,y) \:  = \:  e^{i \eps^{-1} \int_0^t \omega(\eps,s) ds} \Bigl( &
-\eps^{-1} \Bigl( u_s(t,y) \, - \,  u_s(t,a(t)) \, - \, \pa^2_y
u_s(t,a(t)) \frac{y^2}{2} \Bigr) \: \pa_y v_\eps^{sl}(t,y) \\
&  + \: \eps^{-1} \Bigl(
\pa_y u_s(t,a(t)) \,  -  \, \pa^2_y u_s(t,a(t)) y \Bigr) \, v_\eps^{sl}(t,y) \\
& - \, i \, \eps \,  \pa^3_y v_\eps^{reg}(t,y) \: + \: i \,  \pa_t \pa_y\Bigl(v_\eps^{reg}(t,y)
\: + \: v_\eps^{sl}(t,y) \Bigr)  \: + \: O(\eps^{\infty})\Bigr).
\end{align*}
The $O(\eps^{\infty})$ gathers terms with derivatives of
$\varphi$: as the  shear layer profile $V$ decreases exponentially, and the derivatives of
  $\varphi(\cdot-a)$ are supported away from $a$, their  contribution is
indeed exponentially small. Straightforwardly,  
\begin{equation} \label{expo2}
 \|R_\eps(t,\cdot)\|_{W^{0,\infty}_{\alpha}} \: \le  \: C  e^{\frac{\sigma_0 t}{\sqrt{\eps}}},
\end{equation} 
with the same $\sigma_0$ as in \eqref{expo1}.

\bigskip
We are now in a position to prove part i) of Theorem \ref{theorem}. Let us assume {\it
a contrario} that for all $\sigma >0$, 
there exists $m \ge 0,$ $\mu \in [0,1/2)$  and $\delta > 0$   such that 
$$\sup_{0 \le s \le t \le \delta} \| e^{-\sigma (t-s) \sqrt{|\pa_x|}} \, T(t,s) \|_{{\cal L}(H^{m},H^{m-\mu})} \: < \:
  + \infty.$$ 
Let 
$$T_\eps(t,s) : W^{0,\infty}_\alpha(\R_+) \mapsto
W^{0,\infty}_\alpha(\R_+)  $$ 
 the restriction of $T(t,s)$ to the tangential Fourier mode
$\eps^{-1}$. Namely, $T(t,s) \left( e^{i\eps^{-1} x} \,  U_0\right) \: = \: 
e^{i\eps^{-1} x} \, T_\eps(t,s)  U_0$. Similarly, we denote $L_\eps \: = \:
e^{-i\eps^{-1} x} \, L \, e^{i\eps^{-1} x}, \:$ where $L$ is the linearized
Prandtl operator around $u_s$. We have, for all $0 \le s \le t \le \delta$,
$$ \ \| T(t,s) \|_{{\cal L}(W^{0,\infty}_\alpha)}
\: \le \: C \, \eps^{-\mu} \, e^{\frac{\sigma (t-s)}{\sqrt{\eps}}}.$$ 
Let $U = U(t,y)$ the solution of $\pa_t U \: + \: L_\eps U = 0$, that
coincides initially with the  approximation $U_\eps$. On one hand, we get
\begin{equation} \label{upperbound}
  \| U(t,\cdot) \|_{ W^{0,\infty}_\alpha}  \: \le \:  C \, \eps^{-\mu} \,
  e^{\frac{\sigma t}{\sqrt{\eps}}} \, \|
U(0,\cdot) \|_{ W^{0,\infty}_\alpha} \: \le \: C'\,  \eps^{-\mu}  \,
  e^{\frac{\sigma t}{\sqrt{\eps}}}. 
\end{equation}
On the other hand, the difference
$\tilde U = U - U_\eps$ satisfies, for all $t < \delta$, 
$$ \tilde U(t,\cdot) \: = \: \int_0^t T_\eps(t,s) R_\eps(s) ds. $$ 
Estimate \eqref{expo2} implies 
\begin{equation*} 
 \| \tilde U(t,\cdot) \|_{ W^{0,\infty}_\alpha} \: \le \: C \, \eps^{-\mu}
\int_0^t e^{\frac{\sigma (t-s)}{\sqrt{\eps}}} \, e^{\frac{\sigma_0 s}{\sqrt{\eps}}} \, ds  \: \le \: C'  \,
\eps^{1/2-\mu} \,  e^{\frac{\sigma_0 t}{\sqrt{|\eps|}}},
\end{equation*} 
as soon as $\sigma < \sigma_0$. 
Combining this with the estimate \eqref{expo1}, we obtain the lower bound
\begin{align*}
\| U(t,\cdot) \|_{ W^{0,\infty}_\alpha}  \: & \ge \: \| U_\eps(t,\cdot) \|_{
  W^{0,\infty}_\alpha} \: - \:  \| \tilde U(t,\cdot) \|_{W^{0,\infty}_\alpha} \\
& \ge \: c \, e^{\frac{\sigma_0 t}{\sqrt{\eps}}} - C \eps^{\mu-1/2} 
\,  e^{\frac{\sigma_0 t}{\sqrt{\eps}}}
\end{align*}
For $\eps$ small enough, we get
$$ \| U(t,\cdot) \|_{ W^{0,\infty}_\alpha}  \:  \ge \:  c' \,  \,
e^{\frac{\sigma_0 t}{\sqrt{\eps}}}, $$
which contradicts the upper bound \eqref{upperbound},  as soon as 
$\sigma < \sigma_0$ and $t \, \gg \,  \frac{\mu}{\sigma_0 - \sigma} \, |\ln(\eps)| \, \sqrt{\eps}$. 
This achieves the proof of part i).  

\subsection{Stronger ill-posedness for specific $u_s$}
It remains to handle  part ii) of Theorem \ref{theorem}. Roughly, we must find some $u_s$
for which  $e^{-\sigma \sqrt{|\pa_x|}(t-s)} \, T(t,s)$ fails to be  bounded from $H^m$ to
$H^{m-\mu}$,  $\: \mu \ge 0$ arbitrary. Using notations of the previous
paragraph, the keypoint is to build, for any $N$, a growing solution $U_{\eps,N}$ of 
$$\pa_t U_{\eps,N} \: + \: L_\eps U_{\eps,N} \: = \: R_{\eps,N}, \quad \: 
\mbox{ where } \:   \|R_{\eps,N}(t,\cdot)\|_{W^{0,\infty}_{\alpha}} \: \le  \:
C_N \, \left(\eps^N + t^{2N} \right) \, e^{\frac{\sigma_0 t}{\sqrt{\eps}}}. $$
Indeed, we can then  take  $N + 1/2 > \mu$, and conclude along the  same
lines as above. 

\medskip
So far, we have not managed to improve the approximation of the previous
paragraph for general $u_s$. This explains the technical
 restriction $\mu \in [0,1/2)$ of part i). 
In order to obtain a refined approximation, we  consider some special
profiles: we assume that $u_s(0,y) = U_s(y)$, for some exponentially decreasing $U_s$,  satisfying in the
neighborhood of $a > 0$: 
$$  U_s(y) \: = \: U_s''(a) \, \frac{(y-a)^2}{2}, \quad  U_s''(a) < 0.$$ 
Notice that  $a$ is a non-degenerate critical point of $U_s$. 
For such profiles, the approximation   of the previous paragraph reads
$$  U_{\eps}(t,y)  \: = \:  i  \, e^{i \eps^{-1}
  \int_0^t \omega(\eps,s) \,ds}  \, \pa_y \left( v_\eps^{reg}(t,y) \: + \:
v^{sl}_\eps(t,y)\right),  \quad \eps = \frac{1}{n}, \quad  n \in \N_*.$$   
Using that  $U_s$ is quadratic near $y=a$, one can improve this  approximation
through  an expansion of the type
$$ U_{\eps,N}(t,y)  \: = \: U_\eps(t,y) \: + \:  i  \, e^{i \eps^{-1}
  \int_0^t \omega(\eps,s) \,ds}  \: \pa_y \sum_{i=1}^N  \eps^i \, v_\eps^{i,reg}(t,y), $$   
with additional terms $v^{i,reg}_\eps$.
Let us briefly explain the construction of these extra terms. The error
terms  due to $U_\eps$ divide into three categories:
\begin{enumerate}
\item Shear layer terms   involving derivatives of $\varphi$. As mentioned
  before, they are $O\left(\eps^{\infty}\right)$,  and require no correction. 
\item Terms that come from the replacement of $u_s$ by its  Taylor
  expansion in the shear layer equation.They read
\begin{align*}
R_{\eps,1} \: & := \: -\eps^{-1}  e^{i \eps^{-1}
  \int_0^t \omega(\eps,s) \,ds}
\Bigl( u_s(t,y) \, - \,  u_s(t,a(t)) \, - \, \pa^2_y
u_s(t,a(t)) \frac{y^2}{2} \Bigr) \: \pa_y v_\eps^{sl}(t,y), \\
R_{\eps,2} \: & := \: \eps^{-1} e^{i \eps^{-1}
  \int_0^t \omega(\eps,s) \,ds}  \Bigl(
\pa_y u_s(t,a(t)) \,  -  \, \pa^2_y u_s(t,a(t)) y \Bigr) \,
v_\eps^{sl}(t,y). 
\end{align*}
We write 
\begin{align*}
& |R_{\eps,1}| \:  = \: \eps^{-1} e^{\frac{\sigma_0 t}{\sqrt{\eps}}} 
 \left| \int_{a(t)}^y \frac{(z-a(t))^2}{2}
\pa^3_y u_s(t,z) \, dz \right| \: \left| \pa_y v^{sl}_\eps(t,y)\right| \\
&  \le \:  \eps^{-1}  e^{\frac{\sigma_0 t}{\sqrt{\eps}}}  \int_{a(t)}^y \frac{(z-a(t))^2}{2}
  \sum_{k=0}^{2N-1} \frac{t^k}{k!} \, \left|\pa_t^k
\pa_y^{3} u_s(0,z)\right| \, dz \: \left|\pa_y  v_\eps^{sl}(t,y)\right| \: + \:
O(t^{2N})  e^{\frac{\sigma_0 t}{\sqrt{\eps}}} \\
 & \le \:  \eps^{-1} e^{\frac{\sigma_0 t}{\sqrt{\eps}}} \int_{a(t)}^y
\frac{(z-a(t))^2}{2}  \sum_{k=0}^{2N-1}  \frac{t^k}{k!}
\left|\pa_y^{3+2k}U_s(z)\right| \, dz \: \left|\pa_y   v_\eps^{sl}(t,y)\right| \: + \:
O(t^{2N})  e^{\frac{\sigma_0 t}{\sqrt{\eps}}}.
\end{align*}
The second inequality stems from  a Taylor expansion of $u_s$ with respect
to $t$. As $u_s$ satisfies the heat equation, each time derivative can be replaced by two
space derivatives, so the third line. Because  $U_s$ is quadratic in a
vicinity  of $a(t)$ (for short times), and  $v_\eps^{sl}$ and its
derivatives decay exponentially fast, we end up with 
$$ |R_{\eps,1}|  \: \le \: C \, \left( t^{2N} \: + \: \eps^N\right) 
\, e^{\frac{\sigma_0 t}{\sqrt{\eps}}}.$$
A similar bound holds for $R_{\eps,2}$. Hence, these remainders do not
require correction. 
\item Terms that come from the time derivative and the diffusion. We focus
 here on the time derivative, as the diffusion term is simpler and has
 smaller amplitude. This is 
$$ R_{\eps,3} \: := \:   i \, e^{i \eps^{-1}
  \int_0^t \omega(\eps,s) \,ds} \: \pa_t \pa_y\Bigl(v_\eps^{reg}(t,y)
\: + \: v_\eps^{sl}(t,y) \Bigr). $$
Proceeding as for $R_{\eps,2}$, that is with Taylor expansions in $t$,
leads to  
$$   \left| e^{i \eps^{-1}
  \int_0^t \omega(\eps,s) \,ds} \, 
\pa_t \pa_y v_\eps^{sl}(t,y)\right| \: \le \: C \, \left( t^{2N} \: + \:
\eps^N\right)  \, e^{\frac{\sigma_0 t}{\sqrt{\eps}}}.$$
As regards the regular part,
\begin{align*}  
e^{i \eps^{-1}  \int_0^t \omega(\eps,s) \,ds} \: \pa_t \pa_y
v_\eps^{reg}(t,y) \: & = \:  e^{i \eps^{-1}  \int_0^t \omega(\eps,s) \,ds}
\:  H(y-a(t)) \:  \pa_t \pa_y u_s \\
 &  = \:  e^{i \eps^{-1}  \int_0^t \omega(\eps,s) \,ds} \:
 H(y-a(t)) \: F(t,y) \: +
\: O\left(t^{2N} \,   e^{\frac{\sigma_0 t}{\sqrt{\eps}}}\right) 
\end{align*}
where 
$$ F(t,y) \: := \: \sum_{k=0}^{2N-1} \frac{t^k}{k!} \: \pa_y^{3+2k}U_s(y) $$
comes again from a Taylor expansion in $t$. The nice thing about this
$O(1)$  term is that it is identically zero in the vicinity of $y=a(t)$ (for
short times). As a result, the Heaviside function $H(y-a(t))$ in front of
it  does not create any  discontinuity, and  no extra shear layer term is
necessary. One takes care of  this source term by the introduction of 
$$v_\eps^{1,reg}  \: = \: H(y-a(t))\Bigl( u_s(t,y) - \omega(\eps,t) \Bigr)
\int_{a(t)}^y \frac{F(t,z)}{\left(u_s(t,z) - \omega(\eps,t)\right)^2} \, dz,$$  
so that 
$$ U_\eps^1 \: := \: -i \, \eps \, e^{i \eps^{-1}  \int_0^t \omega(\eps,s) \,ds}
\pa_y v_\eps^{1,reg} $$
solves 
$$ \pa_t U_\eps^1 \: + \: L_{\eps} U_\eps^1  \: = \: i \, e^{i \eps^{-1}
  \int_0^t \omega(\eps,s) \,ds} H(y-a(t)) F(t,y) \: + \: O(\eps).$$
 Proceeding recursively, we obtain an approximation as accurate as we
 want.  This ends the proof of the theorem.    
\end{enumerate}

\section{Numerical study} \label{sec4}
In this last  section, we present  numerical illustrations of  the  instability
process. 
\subsection{Numerical test of (SC)}
To check (SC) numerically, it is more convenient to reformulate it with an
{\em Evans function}. We know from Step 2, section \ref{sec2}, that there
are solutions  ${\cal Y}_\pm(z)$  of  \eqref{firstorder} 
satisfying \eqref{ansatz3}-\eqref{ansatz4}. Back to the
ODE \eqref{ODE2}, this yields independent solutions $X_+(\tau,\cdot)$ and 
$X_-(\tau,\cdot)$  respectively growing and decaying as $z$ goes to
$+\infty$. Furthermore, the following asymptotics holds:
\begin{equation*} 
\begin{aligned}
X_\pm(\tau, z) \: & \sim  \: z^{\pm \frac{ i\tau}{2 \lambda}- \frac{7}{2}}
\exp(\pm \frac{1}{2}\lambda z^2), \\    
\pa_z X_\pm(\tau, z) \: & \sim  \: \pm \lambda z^{\pm \frac{i\tau}{2 \lambda}- \frac{3}{2}}
\exp(\pm \frac{1}{2}\lambda z^2),
\end{aligned}
\end{equation*}
with $\lambda = \frac{1-i}{\sqrt{2}}$. Thus, the functions 
$$ W_-(\tau,z) \: := \: \int_z^{+\infty} X_-(\tau,s) \, ds, \quad W_+(\tau,z) \: := \: \int_0^{z}
X_+(\tau,s) \, ds, \quad W_0(\tau,z) = 1, $$
seen as functions of $z$, form a basis of solutions of \eqref{ODE}. They   are respectively
decaying, growing and constant at $+\infty$.   As \eqref{ODE} is preserved by the change of
variable $z \mapsto -z$, the functions $W_-(\tau, -z), \: W_+(\tau, -z), \:
W_0$ form a basis as well. They are respectively decaying, growing and constant at $-\infty$.  
The existence of the heteroclinic orbit is the same as the existence of
some constants $A$ and $B$ such that  
$$ 1 = A W_-(-z,\tau) \: + \:  B W_-(z,\tau)$$ 
for all $z$, or equivalently 
$$ \left( \begin{smallmatrix} 1 \\ 0 \\ 0 \end{smallmatrix} \right) \in
\mbox{ Vect } \left( \left( \begin{smallmatrix} W_-(\tau,0) \\ \pa_z W_-(\tau,0)
  \\ \pa^2_z W_-(\tau, 0) \end{smallmatrix} \right), \,
\left( \begin{smallmatrix} W_-(\tau,0)  \\ -\pa_zW_-(\tau,0)
  \\ \pa^2_z W_-(\tau,0) \end{smallmatrix} \right) \right)$$
This last condition is easily seen to be  equivalent to    
$$  W_-(\tau,0) \: \neq \:  0, \quad \mbox{ and } \pa^2_z W_-(\tau, 0) \: = \: 0.$$
Hence, we must find  $\tau$ with $\Im \tau <
0$ such that   
\begin{equation*}  
\int_0^{+\infty} X_-(\tau,s) \, ds \neq 0, \quad \mbox{ and }  \pa_z
X_-(\tau, 0) \: = \: 0. 
\end{equation*} 
Moreover, we know from Step 3, section \ref{sec2} that for $\Im \tau < 0$, the constraint
$\int_0^{+\infty} X \neq 0$ is satisfied. Finally, the condition  (SC) comes
down to:
$$\pa_z X_-(\tau, 0) \: = \: 0, \quad \mbox{ for some $\tau$ with  $\Im
  \tau < 0$.} $$    

\medskip
To check this, and get a value for $\tau$, one can  use  a shooting method. For any 
  $\tau$  and any $z_0 \gg
 1$,  one can start from the approximation  
 $$ X_-(\tau,z_0) \: \approx \: z_0^{\frac{-i\tau}{2 \lambda}- \frac{7}{2}}
\exp(-\frac{1}{2}\lambda z_0^2), \quad     
\pa_z X_-(\tau, z_0) \: \approx  \: -\lambda z_0^{\frac{-i\tau}{2 \lambda}- \frac{3}{2}}
\exp(-\frac{1}{2}\lambda z_0^2), $$ 
 and integrate backwards \eqref{ODE2} using a Runge-Kutta scheme. 
This gives  access to the function  $\pa_z X_-(\tau,z)$, for any $\tau$ and
any $z \le z_0$. Then,  a Newton-Raphson procedure allows to find a zero in
$\{\Im \tau < 0\}$ for the function $\tau \mapsto  \pa_z
X_-(\tau,0)$. Using such a  procedure, we  have found   
$$\tau \: \approx \:  -0.706 - 0.706 \,i.$$
 Note that this value  is  proportional to  $1+i$, as
expected from the analysis.

\subsection{Simulation of the instability mechanism} 
To observe the instability mechanism described in section \ref{sec1}, we have performed direct
simulations of system \eqref{linear2}. More precisely, we have considered
the velocity
$$ u_s(t,y) \: =  \: U_s(y) \: := \: 2 y \, \exp(-y^2) $$  
 (already studied in \cite{Hon:2003} in the inviscid case), and
 solutions of the type 
$$  u_\eps(t,x,y) \: = \: i \, e^{i \eps^{-1} x} \, \pa_y V_\eps(\eps^{-1} t,y), \quad
 v_\eps(t,y) \: =  \: \eps^{-1} V_\eps(\eps^{-1} t,y). $$ 
The profiles $V_\eps = V_\eps(\theta, y)$ satisfy the singular perturbation problem 
$$    \left(\pa_\theta + i  U_s \right)  \pa_y V_\eps \: - i \: U'_s
V_\eps \: - \:  \eps \, e^{i \eps^{-1} x}  \, \pa^3_y V_\eps \: = \: 0$$
on $V_\eps(\theta,y)$. One more differentiation gives the parabolic like equation
$$   \left(\pa_\theta + i  U_s \right)  \pa^2_y V_\eps \: - i \: U''_s
V_\eps \: - \:  \eps \, \pa^4_y V_\eps \: = \: 0,$$
 fulfilled  with the boundary conditions  
$$ V_\eps\vert_{y=0} \: = \: \pa_y V_\eps\vert_{y=0} \: = \:   \pa_y^3 V_\eps\vert_{y=0} \: = \: 0. $$
We have discretized this equation in space using finite differences on
a stretched grid, and in time
through a Crank-Nicholson scheme. Starting from initial random data
({\it i.e.} with random values at each gridpoint), we have  computed its
time evolution for values of $k=\eps^{-1}$ ranging from $1$ to $3.10^7$. For
sufficiently large times, one  observes that the numerical solution $V^{num}_{\eps}$ behaves like 
$$ V^{num}_\eps(\theta,y) \: \approx \: e^{i \omega^{num}(\eps) \theta} v^{num}_\eps(y)  $$     
in the sense that 
$$ \omega^{num}(\eps) \: := \: \frac{V^{num}_\eps(\theta + \Delta \theta, y) -
  V^{num}_\eps(\theta,y)}{\Delta \theta \,  V^{num}_\eps(t,y)} $$   
gets independent of $\theta$ and $y$. Computations 
show a relation of the type 
$$ \omega^{num}(\eps) \: \sim \: -U_s(a) \: + \: \sqrt{\eps} \, (-0.92 - 0.91 i), $$
see figure~1. Here $a= \frac{1}{\sqrt{2}}$ is as usual the critical point of
$U_s$.  This relation is in very good agreement with the theoretical prediction. 
$$\omega^{th}(\eps) \:  := \:  -U_s(a) \: + \: \eps^{1/2} \,
\frac{| U"_s(a)|^{1/2} \tau}{\sqrt{2}} \: \approx \:   -U_s(a) \: + \: \sqrt{\eps} \, (-0.92 - 0.92 i),$$
if we take for $\tau$ the value $-0.706-0.706i$ found in  the previous
subsection.
\begin{figure} 
\begin{center}
\includegraphics[height=6cm, width=7cm]{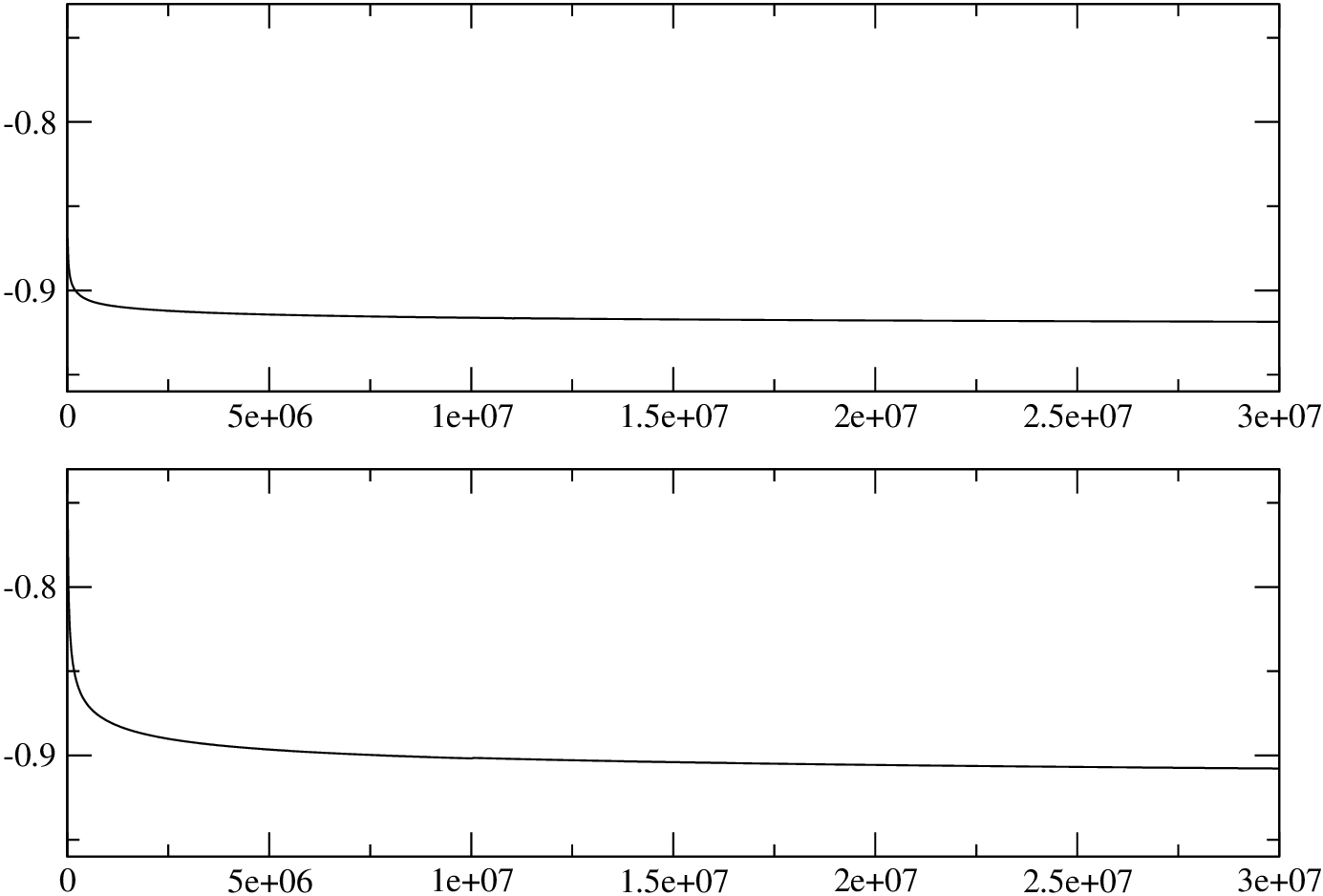}
\caption{The correction to the inviscid eigenvalue. Plot of 
  $\frac{1}{\sqrt{\eps}} \left(\omega^{num}(\eps) + U_s(a)\right)$, seen as a function of
  the tangential frequency $k = \eps^{-1}$. Top: real part. Bottom:
  imaginary part.  As expected from the theory, both approximately
  converge to $-0.9$, as $k$ goes to infinity.}
\end{center}
\label{fig1}
\end{figure}
 Moreover, with this value of $\tau$, one can compute directly the solution
 $V$ of the shear layer equation.  After proper rescaling,  
this  allows for comparison  between the  ``numerical'' and ``theoretical''
eigenmodes.  More precisely, using the notations of \eqref{ansatz2}, one can compare the
functions 
$$ v^{th}_{out}(y) \: := \: \frac{1}{\sqrt{\eps}}
\left(\frac{v^{reg}_{\eps}(y)}{v_\eps^{reg}(\infty)}  - 
 \frac{v_a(y)}{v_a(\infty)} \right)   
\quad \mbox{ and } \quad  v^{num}_{out}(y) \: = \: \frac{1}{\sqrt{\eps}}
\left(\frac{v^{num}_{\eps}(y)}{v_{\eps}^{num}(\infty)}  -  \frac{v_a(y)}{v_a(\infty)} \right)  $$ 
which should both describe the correction to the inviscid eigenmode {\em outside the
shear layer}. As regards the shear layer, one can  compare 
$$ v^{th}_{in}(z) \:  := \:  \frac{1}{\sqrt{\eps}}  \frac{v^{sl}_\eps(\eps^{1/4}(z+a))}{v_\eps^{reg}(\infty)}
 \quad \mbox{ and } \quad  
 v^{num}_{in}(z) \:  := \: \frac{1}{\sqrt{\eps}} 
\frac{v^{num}_{\eps}(\eps^{1/4} (z+a))}{v_{\eps}^{num}(+\infty)}  -
v^{num}_{out}(\eps^{1/4} (z+a)). $$
\begin{figure}
\begin{center}
\includegraphics[height=6cm, width=8cm]{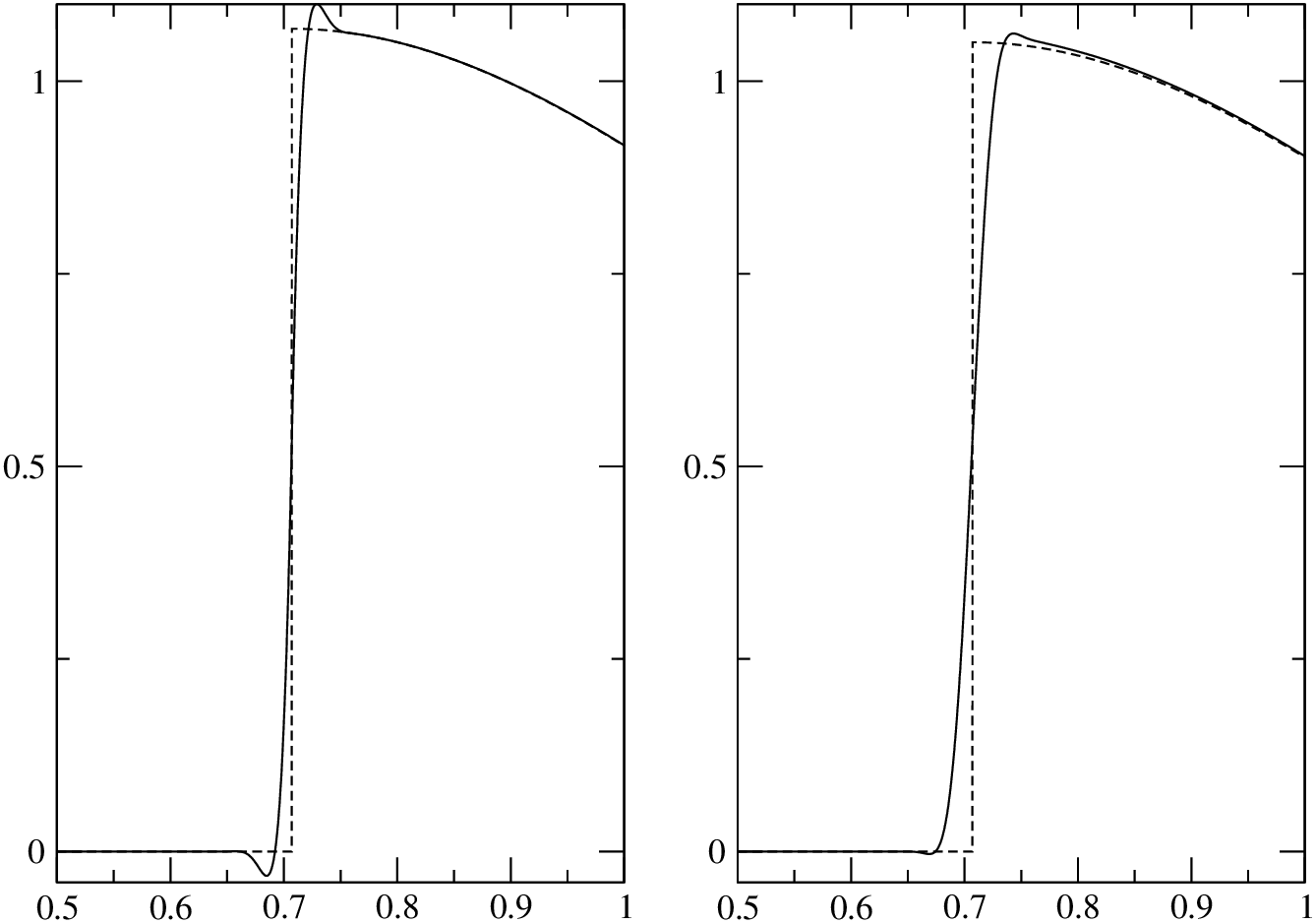}
\caption{Plots of  $v^{th}_{out}$ (dashed line) and $v^{num}_{out}$ (full line),
  seen as functions of $y$, at $\eps = 10^{-7}$. Left and right
  figures correspond respectively to the real and imaginary parts. They
  match, as expected, outside the shear layer.}
\end{center}
\label{fig2}
\end{figure}

\begin{figure}
\begin{center}
\includegraphics[height=6cm, width=8cm]{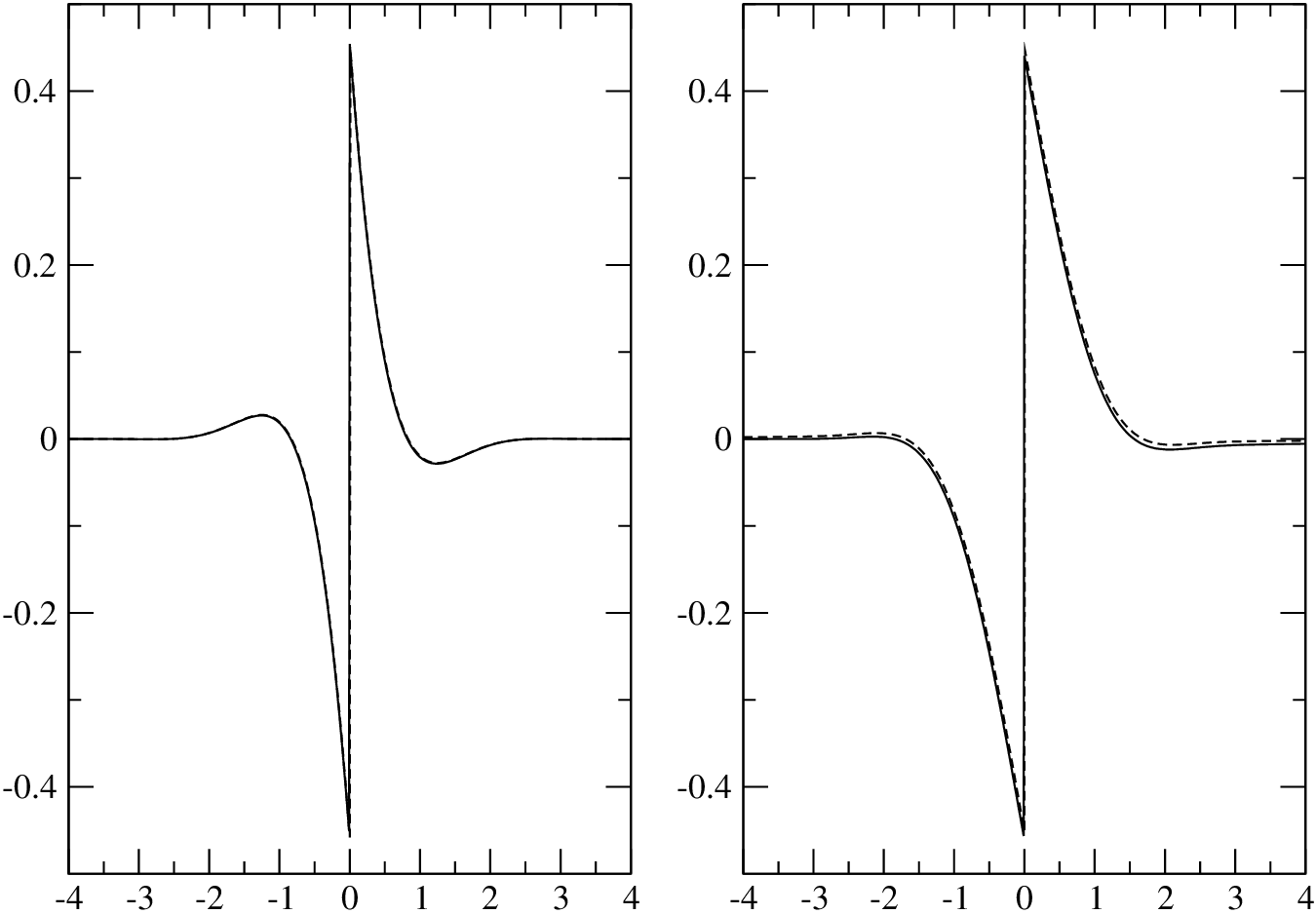}
\caption{
  Plots of  the shear layer corrections $v^{th}_{in}$ (dashed line) and $v^{num}_{in}$ (full line),
  seen as functions of $z$, at $\eps = 10^{-7}$. Left and right
  figures correspond respectively to the real and imaginary parts. }
\end{center}
\label{fig3}
\end{figure}
Illustrations of these comparisons are
given in figures 2 and 3. Again, we obtain an excellent agreement. 
This confirms
that the instability mechanism we have described is indeed effective, and moreover
dominates the  linear dynamics \eqref{linear}.

\section*{Acknowledgements}
The first author would like to thank Guy M\'etivier, Jeffrey Rauch, and
Mikael de la Salle for fruitful discussions. He acknowledges grant ANR-08-JCJC-0104-01.

{\small
\bibliographystyle{acm}
\def\cprime{$'$} \def\cprime{$'$} \def\cprime{$'$}

}

\section*{Appendix : Well-posedness in the analytic setting}
We start from  a simple estimate on the heat equation: 

\medskip
\noindent
{\em  For  $U_0 \in W^{0,\infty}_{\alpha}(\R_+)$,  $\: F \in L^1(0,T; \,
W^{0,\infty}_{\alpha}(\R_+)), $   
the  solution $U$ of 
$$ \pa_t U - \pa^2_y U = F \quad \mbox{on } \: \R \times \R_+,  
 \quad U\vert_{t=0} = U_0, \quad U\vert_{y=0} = 0,$$
satisfies
\begin{equation} \label{estimheat}
 \| U \|_{L^\infty(W_\alpha^{0,\infty})} \: \le \: C \, \left(  \| U_0
\|_{W_\alpha^{0,\infty}} \: + \: \| F \|_{L^1(W_\alpha^{0,\infty})} \right). 
\end{equation}
} 
This estimate follows directly  from the representation formula
$$  U(t,y) \: = \: \int_{\R_+} S(t,y,z) U_0(z) \, dz \: + \: \int_0^t \int_{\R_+} 
S(t-s,y,z) F(s,z) \, dz \, ds $$
 where the heat kernel in the half plane $S(t,y,z)$   reads 
$$ S(t,y,z) \: := \: G(t,y-z) - G(t,y+z), \quad G(t,y) \: := \:
 \frac{1}{\sqrt{4\pi t}}  \exp(-y^2/4t). $$
Details are left to the reader. This estimate  allows to prove proposition
\ref{propanalytic}.
Indeed, by decomposing 
$$u(t,x,y) \: = \: \sum_{k \in \Z} e^{i k x} \, U^k(t,y),$$
 the well-posedness is an easy consequence of  the  {\it a priori} estimate 
$$ \| U^k(t,\cdot) \|_{W^{0,\alpha}} \: \le \: C \, e^{\rho k t}  \, \|
U^k(0,\cdot) \|_{W^{0,\alpha}} $$ 
for some $\rho$. Now, the equation satisfied by $U^k$ is 
$$ \pa_t U^k -  \pa^2_y U^k \: = \: i \, k \, \left( U'_s \int_0^y U^k(t,z)
\, dz \: - \: U_s U^k \right). $$
Using \eqref{estimheat}, we get 
$$ \| U^k(t,\cdot) \|_{W^{0,\infty}_\alpha} \: \le \: C  \| U^k(0,\cdot)
\|_{W^{0,\infty}_\alpha}  \: + \: C_s \, k \, \int_0^t \| U^k(s,\cdot)
\|_{W^{0,\infty}_\alpha} \, ds  $$
where $C_s$ depends on $u_s$. We conclude by the Gronwall lemma.

\end{document}